\documentclass{elsart3-1}

 \usepackage{graphicx}

\usepackage{amssymb}
\usepackage{amsmath}

\usepackage[english,francais]{babel}

\newtheorem{e-proposition}[theorem]{Proposition}

\newtheorem{e-definition}[theorem]{Definition\rm}

\newtheorem{theoreme}{Th\'eor\`eme}[section]

\newtheorem{proposition}[theoreme]{Proposition}

\renewcommand{\theequation}{\arabic{equation}}
\def\og{\leavevmode\raise.3ex\hbox{$\scriptscriptstyle\langle\!\langle$~}}
\def\fg{\leavevmode\raise.3ex\hbox{~$\!\scriptscriptstyle\,\rangle\!\rangle$}}

\journal{the Acad\'emie des sciences}
\begin{document}
\centerline{}
\begin{frontmatter}


\selectlanguage{english}
\title{Multiscale numerical schemes for kinetic equations in the anomalous diffusion limit}


\selectlanguage{english}
\author[crouseilles,ur]{Nicolas Crouseilles},
\ead{nicolas.crouseilles@inria.fr}
\author[ur]{H\'el\`ene Hivert}
\ead{helene.hivert@univ-rennes1.fr}
\author[lemou,ur]{Mohammed Lemou}
\ead{mohammed.lemou@univ-rennes1.fr}

\address[crouseilles]{INRIA}
\address[ur]{IRMAR, Universit\'e de Rennes 1. Campus de Beaulieu. 35000 Rennes.}
\address[lemou]{CNRS}


\medskip
\begin{center}
{\small Received *****; accepted after revision +++++\\
Presented by £££££}
\end{center}

\begin{abstract}
\selectlanguage{english}

We construct  numerical schemes to solve kinetic equations with anomalous diffusion scaling. 
When the equilibrium is heavy-tailed or when the collision frequency degenerates for small velocities, an appropriate 
scaling should be made and the limit model is the so-called anomalous or fractional diffusion model. 
Our first scheme is based on a suitable micro-macro decomposition of the distribution function whereas 
our second scheme relies on a Duhamel formulation of the kinetic equation. Both are \emph{Asymptotic Preserving} (AP): 
they are consistent with the kinetic equation for all fixed value of the scaling parameter $\varepsilon >0$ 
and degenerate into a consistent scheme solving the asymptotic model when $\varepsilon$ tends to $0$.
The second scheme enjoys the stronger property of being uniformly accurate (UA) with respect to $\varepsilon$. 
The usual AP schemes known for the classical diffusion limit cannot be directly applied to the context of anomalous diffusion scaling, 
since they are not able to capture the important effects of large and small velocities. 
We present numerical tests to highlight the efficiency of our schemes. 
{\it To cite this article: N.
Crouseilles, H. Hivert, M. Lemou, C. R. Acad. Sci. Paris, Ser. I 340 (2005) ???????.}
\vskip 0.5\baselineskip

\selectlanguage{francais}

\noindent{\bf R\'esum\'e} \vskip 0.5\baselineskip \noindent
{\bf Sch\'emas num\'eriques multi-\'echelles pour les \'equations cin\'etiques dans la limite de diffusion anormale. }

Nous construisons des sch\'emas num\'eriques pour r\'esoudre les \'equations cin\'etiques dans le r\'egime de diffusion anormale. 
Lorsque l'\'equilibre pr\'esente une queue lourde ou lorsque la fr\'equence de collision d\'eg\'en\`ere pour les petites vitesses, 
un scaling appropri\'e permet d'obtenir un mod\`ele asymptotique appel\'e mod\`ele de diffusion anormale ou fractionnaire. 
Le premier sch\'ema que nous construisons est bas\'e sur une d\'ecomposition micro-macro de la fonction de distribution tandis que 
le second s'appuie sur une formulation de Duhamel de l'\'equation de d\'epart. Ces deux sch\'emas sont \emph{Asymptotic Preserving} (AP) : ils sont consistants avec l'\'equation cin\'etique lorsque le param\`etre d'\'echelle $\varepsilon>0$ est fix\'e et d\'eg\'en\`erent 
en un sch\'ema consistant avec le mod\`ele limite quand $\varepsilon$ tend vers $0$. 
Le deuxi\`eme sch\'ema est m\^eme uniform\'ement pr\'ecis (UA) par rapport \`a $\varepsilon$.
Les sch\'emas AP qui sont connus dans le cas de la limite de diffusion classique ne peuvent pas directement s'appliquer 
au cas de la diffusion anormale car ils ne permettent de capturer les effets importants des petites et des grandes vitesses. 
Nous pr\'esentons des tests num\'eriques pour mettre en \'evidence l'efficacit\'e des sch\'emas que nous pr\'esentons.

{\it Pour citer cet article~: N. Crouseilles, H. Hivert, M. Lemou,  C. R. Acad. Sci.
Paris, Ser. I 340 (2005). ????????????}
\end{abstract}
\end{frontmatter}

\selectlanguage{francais}
\section*{Version fran\c{c}aise abr\'eg\'ee}

Le but de ce travail est de mettre en place des sch\'emas num\'eriques pour les \'equations cin\'etiques lin\'eaires dans le cas de la limite de diffusion anomale. Quand la distribution d'\'equilibre est une fonction \`a queue lourde ou lorsque la fr\'equence de collision est d\'eg\'en\'er\'ee 
pour les petites vitesses, l'\'equation cin\'etique \eqref{KinEq} tend vers l'\'equation dite de diffusion anormale \eqref{LimEq} 
lorsque le param\`etre d'\'echelle $\varepsilon$ tend vers z\'ero. 
Dans cette limite, une raideur appara\^it dans l'\'equation de d\'epart 
et une r\'esolution num\'erique directe du probl\`eme peut devenir extr\^emement co\^uteuse puisque les param\`etres num\'eriques 
doivent a priori \^etre adapt\'es \`a $\varepsilon$.  
La construction de sch\'emas dits \emph{Asymptotic Preserving} (AP) permet de r\'epondre \`a cette contrainte : 
ces sch\'emas restent consistants avec l'\'equation cin\'etique tout en s'affranchissant de la contrainte sur les param\`etres 
num\'eriques, et d\'eg\'en\`erent vers l'\'equation limite quand $\varepsilon$ tend vers $0$. 

Dans le cas de l'asymptotique de diffusion anormale, la mise en place de tels sch\'emas s'av\`ere plus compliqu\'ee que dans le cas classique. 
En effet, en plus de la raideur \'evoqu\'ee ci-dessus ($\varepsilon$ tend vers $0$), il est crucial de capturer les effets des grandes 
et petites vitesses pour que ces sch\'emas d\'eg\'en\`erent vers des approximations consistantes de l'\'equation de diffusion anormale \eqref{LimEq}. 
Un sch\'ema num\'erique inspir\'e des approches AP standards ne tiendrait pas compte de ces effets, et d\'eg\'en\`ererait  
vers une approximation d'une \'equation de diffusion classique et non vers celle du mod\`ele correct de diffusion anormale. 

Deux cas sont consid\'er\'es dans ce travail : le cas d'un \'equilibre \`a queue lourde et celui d'une fr\'equence de collision d\'eg\'en\'er\'ee en $0$. 
Dans les deux cas, nous construisons d'abord un sch\'ema bas\'e sur une 
d\'ecomposition micro-macro de la solution, qui fournit un sch\'ema multi-\'echelle pour \eqref{KinEq} compl\`etement explicite en temps ; 
ensuite, nous pr\'esentons un sch\'ema  
bas\'e sur une formulation de Duhamel de l'\'equation cin\'etique \eqref{KinEq} ayant une propri\'et\'e plus forte que la propri\'et\'e AP : 
la  pr\'ecision de ce sch\'ema est uniforme (UA) par rapport \`a $\varepsilon$. 
 
 
Dans chaque cas, l'\'ecriture d'un sch\'ema semi-discret en temps permet d'obtenir une formulation qui tend vers 
l'\'equation de diffusion anormale lorsque $\varepsilon$ tend vers $0$ si l'espace des vitesses est consid\'er\'e comme continu et les int\'egrations en vitesses sont r\'ealis\'ees exactement. Cependant, une discr\'etisation directe de l'espace des vitesses pour r\'ealiser les int\'egrations num\'eriquement 
ne permet pas de prendre en compte les effets des grandes et petites vitesses, qui sont \`a  l'origine de la limite de diffusion anormale. 
Dans ce travail, nous montrons qu'il est donc n\'ecessaire d'effectuer des transformations sur certaines int\'egrales en vitesses avant de les discr\'etiser. 
En l'occurrence, nous effectuons des changements de variables ad\'equats pour faire appara\^itre naturellement les termes \`a l'origine de l'\'equation asymptotique. 
Nous obtenons ainsi des sch\'emas compl\`etement discr\'etis\'es ayant la propri\'et\'e AP et UA. 
Nous pr\'esentons \'egalement des tests num\'eriques qui mettent en \'evidence la limite de diffusion anormale de nos sch\'emas quand $\varepsilon$ tend vers $0$.
Cette note est une version abr\'eg\'ee de \cite{CrousHivertLemou1}, \cite{CrousHivertLemou2}.



\selectlanguage{english}
\section{Introduction}
We consider the kinetic equation
\begin{equation}
\label{KinEq}
\varepsilon^\alpha\partial_t f +\varepsilon v\cdot \nabla_x f =Q(f) := \nu(v)\left(\rho_\nu M-f\right), \;\;\;\;\;\; f(0,x,v)=f_0(x,v),  
\end{equation}
where $f$ is a distribution function which depends on the time $t\ge 0$, the space variable $x\in  \mathbb{R}^d$ and the velocity $v\in \mathbb{R}^d$ with $d=1,2,3$, $f_0$ is a given initial data. In the sequel, we will denote by brackets the integration in $v$. We define the density $\rho$ by 
$\rho(t,x)=\left\langle f(t,x,v)\right\rangle =: \int_{v\in\mathbb{R}^d} f(t,x,v)\mathrm{d}v$, and the quantity $\rho_\nu$ 
(which is not the usual density) by $\rho_\nu(t, x)=\left\langle \nu(v) f\right\rangle/\left\langle \nu(v)M(v)\right\rangle$. Note that 
this definition of $\rho_\nu$ ensures the local mass conservation $\langle Q(f)\rangle =0$. 
The positive number $\varepsilon$ is a scaling parameter and $\alpha$ is a power which will be chosen 
according to the physical nature of the problem (see below for more details).  
The equilibrium distribution function $M$  is a normalized even function of $v$. 
 We will consider two physically relevant cases, depending on the nature of the equilibrium function $M$ (see \cite{MelletMischlerMouhot}, \cite{Puel2}) 
 and of the collision frequency $\nu$ (see \cite{Puel1}): 
\begin{enumerate}
\item \textbf{Case $1$ - heavy-tail:}  $M$ is a power-tailed equilibrium and $\nu(v)=1$. To simplify, we will consider the case $M(v)=m/(1+|v|^\beta)$ 
where $\beta\in(d,d+2)$ and the normalization parameter $m$ is chosen such that $\left\langle M\right\rangle = 1$. In this case, the appropriate 
choice of $\alpha$ is:  $\alpha=\beta-d\in (0, 2)$. 
\item \textbf{Case $2$ - degenerate collision frequency:}  $\nu(v)\underset{v\to 0}{\sim} \nu_0 |v|^{d+2+\beta}, \beta>0, \nu_0>0$ and 
$M(v)=\frac{\mathrm{e}^{-|v|^2/2}}{\sqrt{2\pi}}$. To simplify, we will consider the case $\nu(v)=\nu_0 |v|^{d+2+\beta}$ for some $\beta>0$. 
In this case, the appropriate 
choice of $\alpha$ is: $\alpha=(2+2d+\beta)/(1+d+\beta) \in (1,2)$.
\end{enumerate}
Note that in these two cases, the quantity $D=\left\langle |v|^2 \frac{M}{\nu(v)}\right\rangle$ which usually appears as the diffusion coefficient in the classical 
diffusion case, is not finite.  In fact, in our context, the limit for small $\varepsilon$ of  \eqref{KinEq} is not given by diffusion but by \emph{anomalous diffusion} 
equation, which can be written with a fractional Laplacian 
\begin{equation}
\label{LimEq}
\partial_t \hat{\rho}+\kappa |k|^\alpha \hat{\rho}=0, ~~\hat{\rho}(0,k)=\left\langle \hat{f}_0(k,v)\right\rangle. 
\end{equation}
In \eqref{LimEq}, the quantity $\hat{\rho}$ (respectively $\hat{f}$) denotes the space Fourier transform of the function 
$\rho$ (respectively $f$) and $k$ is the Fourier variable. 
The coefficient $\kappa$ can be expressed in the two cases: in Case $1$  of heavy-tail equilibrium, it writes
\[
\kappa=\left\langle \frac{m}{|v|^\beta}\frac{(v\cdot e)^2}{1+(v\cdot e)^2}\right \rangle,
\]
and in Case $2$ of the degenerate collision frequency, it writes  
\[
\kappa=\frac{1}{\sqrt{2\pi}}\frac{~\nu_0^{1-\alpha}}{d+1+\beta}\left\langle \frac{1}{|v|^\beta}\frac{(v\cdot e)^2}{1+(v\cdot e)^2}\right \rangle,
\]
where $e$ denotes any unitary vector of $\mathbb{R}^d$. 
The anomalous diffusion limit for kinetic equations has been studied in \cite{Puel1} 
in the case of degenerate collision frequency and in \cite{Puel2}, \cite{MelletMischlerMouhot} in the case of heavy-tailed equilibrium. 

\section{Micro-macro scheme}

In this section, we derive an AP micro-macro scheme for \eqref{KinEq} in the case of the anomalous diffusion limit.  
This approach bears similarities with the one developed in \cite{LemouMieussens} and has the strong advantage 
of treating the transport explicitly. Note that an implicit AP scheme can be constructed directly, 
with no use of the micro-macro decomposition. Indeed we proceed as follows. We first start by an implicit formulation of \eqref{KinEq}
\begin{equation}
\label{implicit}
f^{n+1}=\lambda(v)\left( I+\frac{\varepsilon \lambda(v)}{\nu(v)}v\cdot \nabla_x \right)^{-1} \rho^{n+1}_\nu M +(1-\lambda(v))\left( I+\frac{\varepsilon \lambda(v)}{\nu(v)}v\cdot \nabla_x \right)^{-1}f^n,
\end{equation}
where $\lambda(v)=\Delta t \nu(v) /(\varepsilon^\alpha+\Delta t\nu(v))$,  $\Delta t$ is the time step, 
$f^n\sim f(t_n)$, $\rho^n=\left\langle f^n\right\rangle$ and $g^n=f^n-\rho^n M$ such that $\left\langle g^n\right\rangle=0$. 
We then integrate \eqref{implicit} against $\nu$ with respect to $v$ and get an expression of $\rho_\nu^{n+1}$; then, we use  an adequate treatment of the integrals in velocity appearing in the so-obtained expression of $\rho_\nu^{n+1}$ to ensure the AP property of the scheme, 
using adapted changes of variables (see \cite{CrousHivertLemou1} for more details). 
Once this expression of $\rho_\nu^{n+1}$ is obtained, we report it in \eqref{implicit} to compute $f^{n+1}$. 
Here, we are interested in deriving numerical schemes where the transport part is explicit in time. To this end, we use a suitable micro-macro decomposition which leads to the following numerical scheme (see again \cite{CrousHivertLemou1} for more details). 

\begin{proposition}
We introduce the semi-discrete micro-macro scheme defined for all $x\in\mathbb{R}^d, v\in \mathbb{R}^d$ and all time index $0\le n\le N$, with $N\Delta t=T$ 
($T>0$), by 
\begin{align}
\label{mm1}
\frac{\rho^{n+1}-\rho^n}{\Delta t}&+\mathcal{F}^{-1} \left(\mathcal{I} \hat{\rho}^{n+1/2}\right) +\mathcal{F}^{-1} \left( \mathcal{I} \frac{\left\langle \nu(v) \hat{g}^{n+1}\right\rangle}{\left\langle \nu(v) M\right\rangle} \right)
+ \varepsilon \left\langle \frac{v\cdot \nabla_x g^n}{\varepsilon^\alpha + \nu(v)\Delta t}\right\rangle =0, \\
\label{mm2}
\displaystyle\frac{g^{n+1}-g^n}{\Delta t}\displaystyle
&+\varepsilon^{1-\alpha}v\cdot \nabla_x \rho^n M+\varepsilon^{1-\alpha}\left( v\cdot \nabla_x g^n-\left\langle v\cdot \nabla_x g^n \right\rangle M \right)
=-\frac{\nu(v)}{\varepsilon^\alpha}\left(g^{n+1}-\frac{\left\langle\nu(v)g^{n+1}\right\rangle}{\left\langle \nu(v)M\right\rangle}M\right),
\end{align}
where $\mathcal{F}^{-1}$ denotes the inverse Fourier transform in space. The quantity 
$\hat{\rho}^{n+1/2}$ can be chosen equal to $\hat{\rho}^n$ or  to $\hat{\rho}^{n+1}$ depending on the desired asymptotic scheme (explicit   or implicit in time) and  the quantity $\mathcal{I}$  is given by
\begin{equation}
\label{I}
\mathcal{I}= \varepsilon^{1-\alpha}\left\langle \lambda(v) \frac{\mathrm{i} k\cdot v }{1+\frac{\varepsilon \lambda(v)}{\nu(v)}\mathrm{i} k\cdot v} M\right\rangle 
=\varepsilon^{2-\alpha}\left\langle  \frac{ \nu(v)  \lambda(v)^2 (k\cdot v)^2}{\nu(v)^2+\varepsilon^2 \lambda(v)^2 (k\cdot v)^2} M\right\rangle.  
\end{equation}
This scheme is of order $1$ for any fixed $\varepsilon>0$ and enjoys the AP property: for a fixed $\Delta t$, the scheme 
degenerates into a first order in time scheme for \eqref{LimEq} when $\varepsilon$ goes to zero. 
\end{proposition}

Let us say a few words about the derivation of \eqref{mm1}-\eqref{mm2}. 
Since $\lambda=\mathcal{O}(\Delta t)$, we have from \eqref{implicit} 
\begin{equation}
\label{implicit2}
f^{n+1}=\lambda(v)\left( I+\frac{\varepsilon \lambda(v)}{\nu(v)}v\cdot \nabla_x \right)^{-1} \rho^{n+1}_\nu M 
+(1-\lambda(v))f^n + \mathcal{O}(\Delta t). 
\end{equation}
Then, we integrate \eqref{KinEq} in $v$ to get the continuity equation on $\rho$ 
and write a Euler implicit scheme: $\rho^{n+1} = \rho^n - \Delta t \varepsilon^{1-\alpha} \nabla_x\cdot\langle v f^{n+1}\rangle$.  
Then, we replace $f^{n+1}$ by \eqref{implicit2} in this scheme, use the identity 
$\rho_\nu^{n+1} = \rho^{n+1} + \langle \nu g^{n+1}\rangle /\langle \nu M \rangle$ 
and the evenness of $M(v)$ to get \eqref{mm1}. To get \eqref{mm2}, we just apply $(I-\Pi)$ to \eqref{KinEq} 
where $\Pi f=\langle f\rangle M$ and discretize the obtained equation on $g$ by an explicit scheme for the transport terms 
and an implicit scheme for the collision terms.

Before discussing the delicate issue relative to the velocity discretization in  \eqref{mm1}-\eqref{mm2}, we first briefly explain how $(\rho^{n}, g^n)$ 
are computed recursively from \eqref{mm1}-\eqref{mm2}, assuming that the space and velocity discretizations have been already fixed. 
The idea is to start with  \eqref{mm2} to find an expression for $g^{n+1}$. In Case $1$ (heavy-tailed equilibrium) as $\nu(v)=1$, the term $\left\langle \nu(v)g^{n+1}\right\rangle$ in the right hand side of \eqref{mm2} vanishes, giving easily an expression for $g^{n+1}$, which is then reported in \eqref{mm1} and so on. 
However, in Case $2$ (degenerate collision frequency), it is necessary to extract $\left\langle \nu(v) g^{n+1}\right\rangle$ 
before solving the equation in $g^{n+1}$. To do that, we express $g^{n+1}$ from \eqref{mm2} in terms of $\rho^n, g^n$ and $\langle \nu(v) g^{n+1}\rangle$; 
we multiply this obtained expression of $g^{n+1}$ by $\nu(v)$ and  integrate in $v$ to get the following expression for $\left\langle \nu(v) g^{n+1}\right\rangle$ 
\[
\left\langle \nu(v)g^{n+1}\right\rangle=
\frac{ \displaystyle
  \left\langle 
    \frac{\nu(v) g^n}{1+\frac{\Delta t}{\varepsilon^\alpha}\nu(v)}
  \right\rangle 
- \Delta t 
  \left\langle 
    \frac{\nu(v) \left( \frac{\varepsilon}{\varepsilon^\alpha}v\cdot \nabla_x \rho^n M+\frac{\varepsilon}{\varepsilon^\alpha}\left( v\cdot \nabla_x g^n-\left\langle v\cdot \nabla_x g^n \right\rangle M \right) \right)}{1+\frac{\Delta t}{\varepsilon^\alpha} \nu(v) }
  \right\rangle 
}{ \displaystyle
1-\frac{\Delta t }{\varepsilon^\alpha \left\langle \nu(v)M \right\rangle}\left\langle  \frac{\nu(v)^2M }{1+\frac{\Delta t}{\varepsilon^\alpha}\nu(v)} \right\rangle
}.
\]
Reporting this in \eqref{mm2} enables to compute $g^{n+1}$, which is reported into \eqref{mm1} to get $\rho^{n+1}$.

Now, we come to the construction of complete discretization from \eqref{mm1}-\eqref{mm2} and in particular the important point of 
discrete integrations in velocity. Indeed, the behavior of the scheme when $\varepsilon$ goes to $0$ is intimately related to 
the discretization we make to approximate \eqref{I}. To see this, we observe that a direct discretization (by rectangle formula for instance) 
of $\mathcal{I}$ given by \eqref{I} converges to $0$ when $\varepsilon$ goes to $0$ since the numerical 
approximation of the bracket in  \eqref{I} is always finite. Therefore this leads to the wrong limit and to overcome this problem, we proceed as follows. 
The general idea is to perform a suitable change of variables in \eqref{I} before discretizing it in velocity. 
In Case $1$ (heavy-tailed equilibrium) we make the change of variable  $w=\varepsilon\lambda |k| v$ in $\mathcal{I}$ before applying the discretization and in Case $2$ (degenerate collision frequency) we make the change of variables $w=\varepsilon |k| v/\nu(v)$ in $\mathcal{I}$  before discretizing the brackets. Hence, when $\varepsilon$ tends to $0$, $\mathcal{I}$ degenerates into the discretized coefficient $\kappa$ of \eqref{LimEq}. This ensures the AP property of the fully discretized scheme. 
Note that the other 
velocity integrations in \eqref{mm1}-\eqref{mm2} are discretized directly without any change of variables. 
Finally, a standard upwind scheme is used for the spatial discretization of \eqref{mm2}.  

\begin{figure}[!htbp]
\begin{center}
\begin{tabular}{@{}c@{}c@{}}
\includegraphics[width=6.5cm]{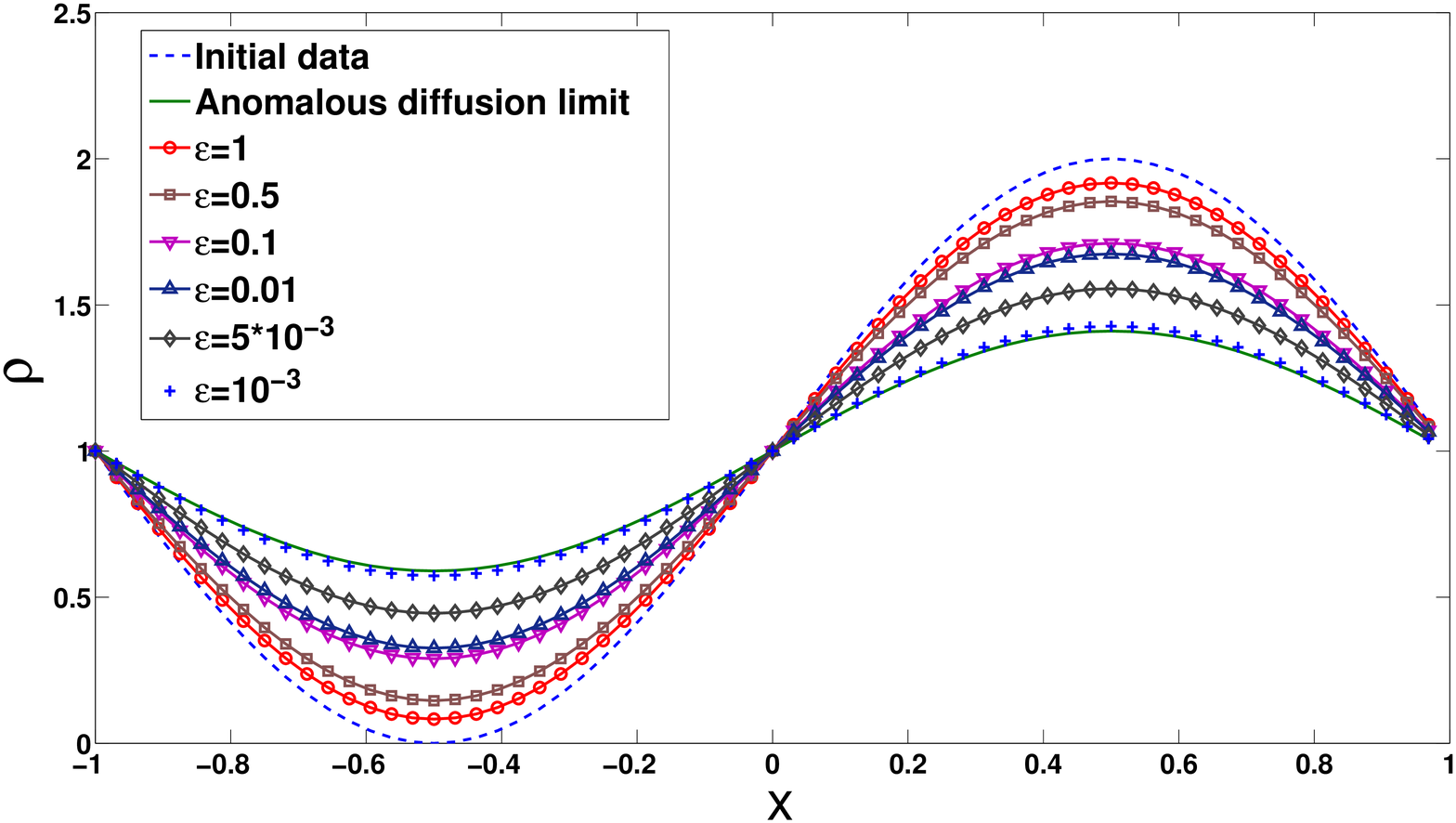} &     
\includegraphics[width=6.5cm]{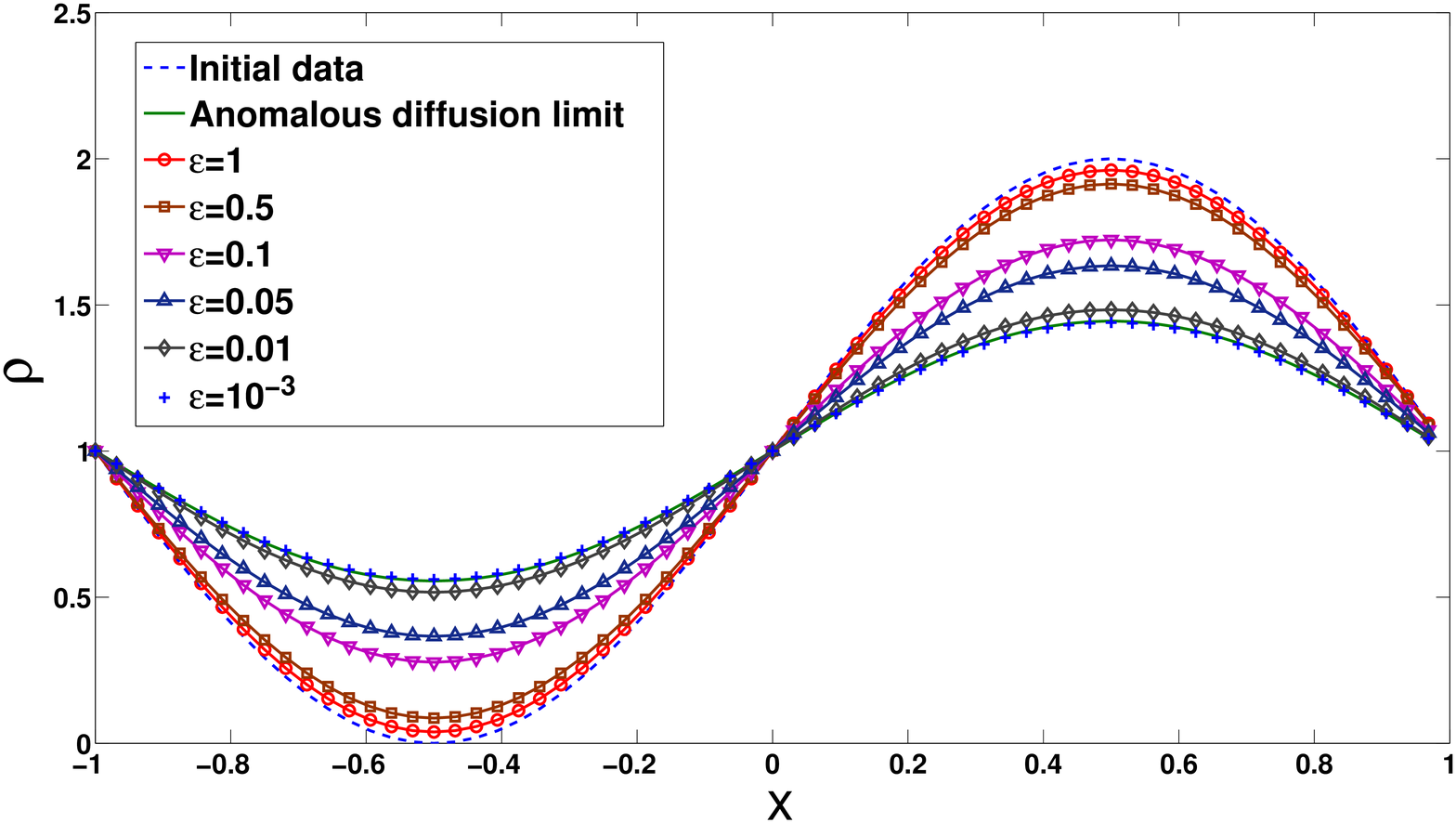}
\end{tabular}
\caption{The densities $\rho(T=0.1, x)$ for the micro-macro schemes and the anomalous diffusion limit scheme. Left: the case of heavy-tailed equilibrium. Right: the case of degenerate collision frequency.}
\label{figure1}
\end{center}
\end{figure}

To highlight the AP character of this scheme, we present in Fig. \ref{figure1} the densities $\rho$ computed at time $T=0.1$ with this scheme for some $\varepsilon$. 
We took the initial data $f_0(x,v)=(1+\sin(\pi x))M(v)$ for $x\in[-1,1]$ and $64$ points of discretization in space. 
In Case $1$, we took $\beta=2.5$ and in Case $2$ we took $\beta=0.5$. In both cases we considered $v_{\text{max}}=5$, and $200$ points of discretization in velocity and $\Delta t=10^{-3}$.

\section{Duhamel formulation based scheme}

In this section, our goal is to derive numerical schemes for \eqref{KinEq} which go beyond the AP property. 
More precisely, our aim to construct numerical schemes whose accuracy in time is uniform with respect to $\varepsilon$. 
Note that although the scheme of the last section is AP, it does not enjoy this uniform accuracy property (see \cite{CrousHivertLemou1}). 
The starting point of the derivation is a Duhamel formulation of \eqref{KinEq} in Fourier variable, which leads,  after 
an integration in $v$ against $\nu$, to an exact equation on $\hat{\rho}_\nu$ 
\begin{equation}
\label{Duhamel_formula}
\hat{\rho}_\nu(t,k)=\hat{A}_0(t,k) +(\left\langle \nu(v)M\right\rangle)^{-1}\int_0^{\frac{t}{\varepsilon^\alpha}}\left\langle \mathrm{e}^{-s(\nu(v)+\mathrm{i} \varepsilon  k\cdot v)}\nu(v)^2M(v)\right\rangle \hat{\rho}_\nu(t-\varepsilon^\alpha s,k)\mathrm{d} s, 
\end{equation}
with $\hat{A}_0(t,k)=\left\langle \mathrm{e}^{-\frac{t}{\varepsilon^\alpha}\left( \nu(v)+\mathrm{i} \varepsilon k\cdot v \right)}\nu(v) \hat{f}_0(k,v)\right\rangle/\left\langle \nu(v)M\right\rangle$.
Considering discrete time values $t_n=n\Delta t, 0\le n\le N, N\Delta t=T, (T,\Delta t>0)$, and evaluating \eqref{Duhamel_formula} 
at $t=t_{n+1}$ leads to 
$$
\hat{\rho}_\nu(t_{n+1},k)=\hat{A}_0(t_{n+1},k) +(\left\langle \nu(v)M\right\rangle)^{-1}\sum_{j=0}^n\int_{\frac{t_{j}}{\varepsilon^\alpha}}^{\frac{t_{j+1}}{\varepsilon^\alpha}}\left\langle \mathrm{e}^{-s(\nu(v)+\mathrm{i} \varepsilon  k\cdot v)}\nu(v)^2M(v)\right\rangle \hat{\rho}_\nu(t_{n+1}-\varepsilon^\alpha s,k)\mathrm{d} s. 
$$
For $s\in [t_j/\varepsilon^\alpha, t_{j+1}/\varepsilon^\alpha]$, using $\hat{\rho}_\nu^n\sim\hat{\rho}_\nu(t_n)$, we consider the approximation 
$\rho_\nu(t_{n+1}-\varepsilon^\alpha s,k) \approx (t_{j+1}-\varepsilon^\alpha s)/\Delta t \; \rho_\nu^{n+1-j} + (\varepsilon^\alpha s-t_j)/\Delta t \; \rho_\nu^{n-j}$ and get the following proposition. 
\begin{proposition}
We introduce the Duhamel formulation based scheme defined for all $x\in \mathbb{R}^d, v\in\mathbb{R}^d$ and all time index $0\le n\le N, N\Delta t=T$ by 
\begin{equation}
\label{Duhamel_scheme}
\hat{\rho}^{n+1}_\nu=\hat{A}_0(t^{n+1},k)+\frac{1}{\left\langle \nu(v)M\right\rangle}\sum\limits_{j=0}^n (c_j-b_j)\hat{\rho}^{n+1-j}_{\nu}+b_j\hat{\rho}^{n-j}_\nu,
\end{equation}
where
\begin{align}
b_j&=\int_\frac{t_j}{\varepsilon^\alpha}^\frac{t_{j+1}}{\varepsilon^\alpha} \frac{\varepsilon^\alpha s-t_j}{\Delta t}\left \langle \nu(v)^2M \mathrm{e}^{-s\nu(v)} \right\rangle \mathrm{d}s
 +\int_\frac{t_j}{\varepsilon^\alpha}^\frac{t_{j+1}}{\varepsilon^\alpha} \frac{\varepsilon^\alpha s-t_j}{\Delta t}\left \langle \nu(v)^2 M \left( \mathrm{e}^{-s(\nu(v)+\mathrm{i}\varepsilon k\cdot v)}- \mathrm{e}^{-s\nu(v)} \right)\right\rangle \mathrm{d}s, \label{bj} \\
 c_j&=\int_\frac{t_j}{\varepsilon^\alpha}^\frac{t_{j+1}}{\varepsilon^\alpha} \left \langle \nu(v)^2 M \mathrm{e}^{-s\nu(v)} \right\rangle \mathrm{d}s
 +\int_\frac{t_j}{\varepsilon^\alpha}^\frac{t_{j+1}}{\varepsilon^\alpha}\left \langle \nu(v)^2 M \left( \mathrm{e}^{-s(\nu(v)+\mathrm{i}\varepsilon k\cdot v)}- \mathrm{e}^{-s\nu(v)} \right)\right\rangle \mathrm{d}s. \label{cj}
\end{align}
This scheme is of order $1$ uniformly in $\varepsilon$: $\exists C>0$ independent from $\varepsilon$, such that 
$\max_{k} |\hat{\rho}_\nu(t^n, k)-\hat{\rho}_\nu^n(k)| \leq C \Delta t, \forall 0\leq n\leq N, N\Delta t=T$.  
\end{proposition}
We refer to \cite{CrousHivertLemou1} and \cite{CrousHivertLemou2} for more details and for the proof of the uniform accuracy of this scheme.

The result of this proposition implies in particular, that this semi-discrete scheme enjoys the AP property. 
Moreover, we will show below how to discretize the velocity integrals in order to preserve this uniform accuracy property 
for the fully discretized scheme. In fact, as in the previous section, a direct rectangular velocity discretization in \eqref{bj}-\eqref{cj} 
leads to a wrong limit as $\varepsilon$ goes to $0$. Therefore, a change of variable in these velocity integrals is necessary 
before any discretization 
to ensure the right anomalous diffusion limit. 

In Case $1$ (heavy-tailed equilibrium), we perform the change of variables $w=\varepsilon|k|v$ 
only in the second velocity bracket of \eqref{bj} and do the same for \eqref{cj}. Then, a rectangle formula 
is applied to the obtained integrals. The first velocity bracket in \eqref{bj} 
is computed directly using a rectangle formula, and the same is done for \eqref{cj}. 
We proceed similarly for Case $2$ (degenerate collision frequency), except that we perform the change of variable $w=\varepsilon |k| v/\nu(v)$ 
at the same places.  
It is important to note that the discretizations in velocity are independent of $\varepsilon$ and 
we still have the first order (in time) uniform accuracy with respect to $\varepsilon$; 
this statement is proved in \cite{CrousHivertLemou2}. 

\begin{figure}[!htbp]
\begin{center}
\begin{tabular}{@{}c@{}c@{}}
\includegraphics[width=6.5cm]{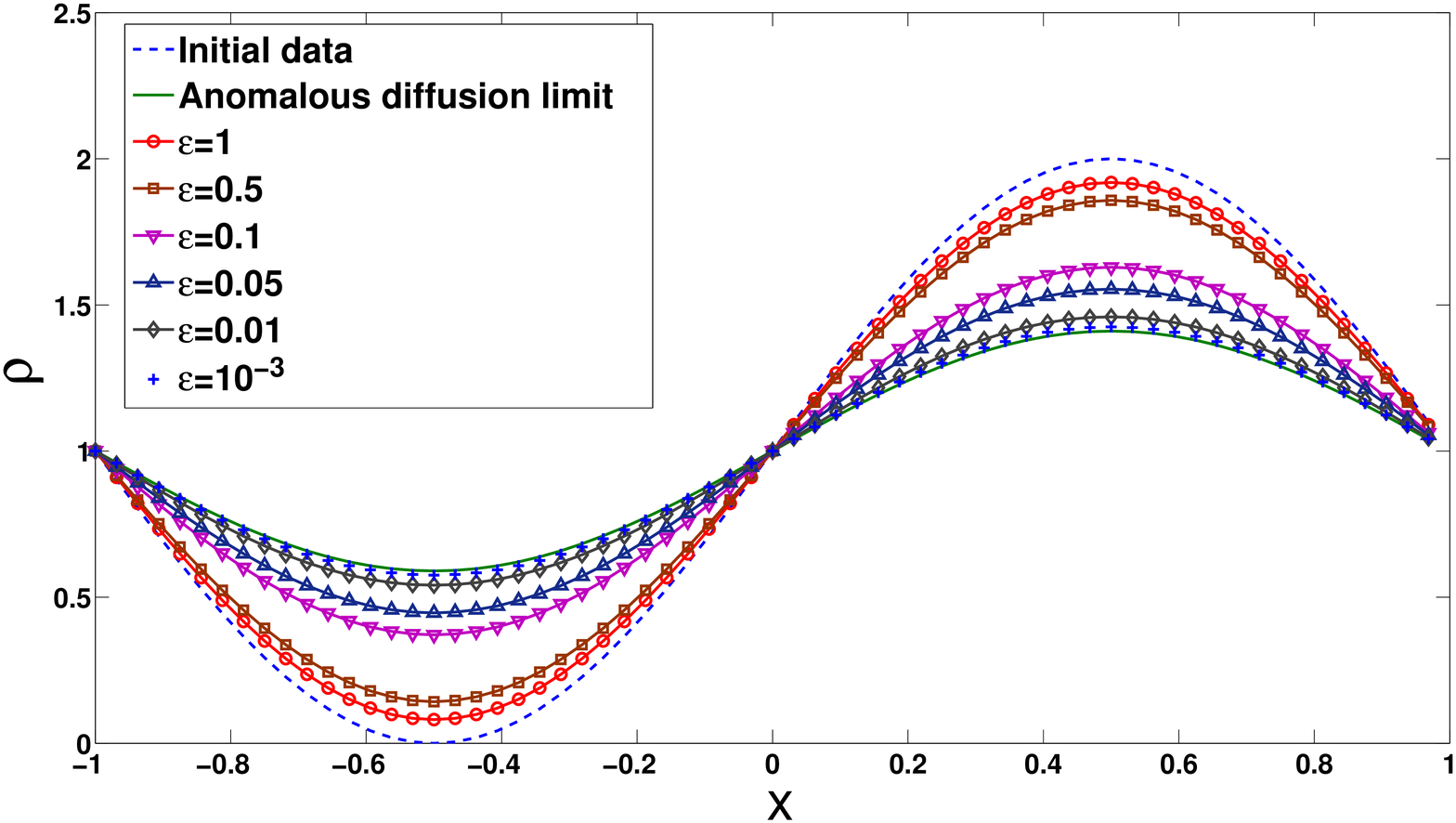} &     
\includegraphics[width=6.5cm]{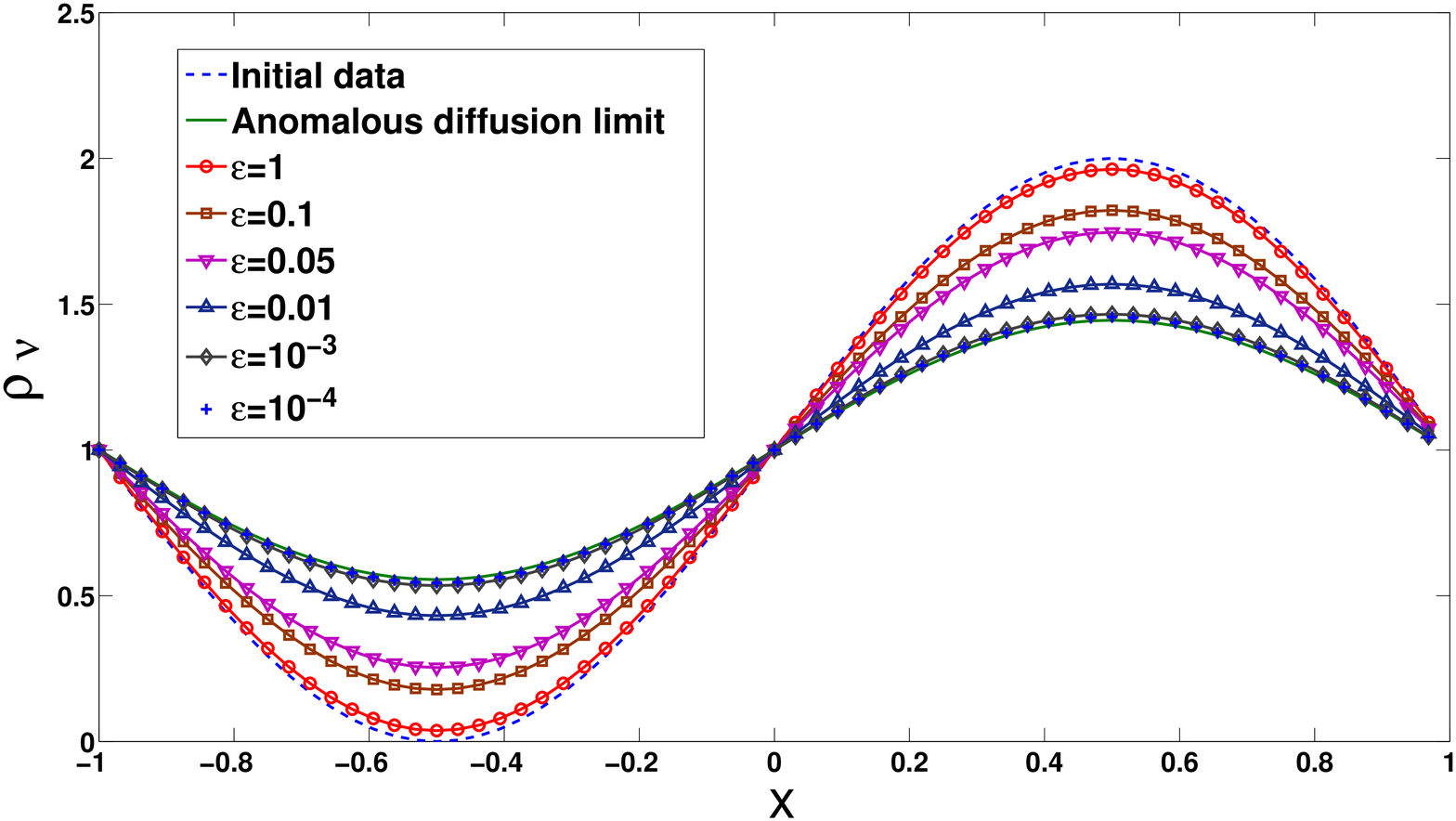}
\end{tabular}
\caption{The densities $\rho_\nu(T=0.1, x)$ for the Duhamel schemes and the anomalous diffusion limit scheme. Left: the case of heavy-tailed equilibrium. Right: the case of degenerate collision frequency.}
\label{figure2}
\end{center}
\end{figure}

In order to highlight the AP character of this scheme, we present in Fig. \ref{figure2} 
the densities obtained  with the initial data $f_0(x,v)=(1+\sin(\pi x))M(v), x\in[-1,1]$ at time 
$T=0.1$ and with $\Delta t=10^{-2}$. In the case of the heavy-tailed equilibrium, we took $\beta=2.5$ and in the case of the degenerate collision frequency, we took $\beta=0.5$. In both cases, we consider $v_{\text{max}}=5$ and $200$ points of discretization in velocity.

Note that in 
both cases, the distribution function $f$ 
can be easily recovered from the values of $\rho_\nu$ given by the previous scheme,
using a Duhamel formulation of the kinetic equation for $f$. The same uniform accuracy property for $f$ 
is inherited for that of $\rho_\nu$. 

\section*{Acknowledgements}
We acknowledge support by the ANR project Moonrise (ANR-14-CE23-0007-01). 
This work was also partly supported by the ERC Starting Grant Project GEOPARDI.


\begin{thebibliography}{00}

 \bibitem{Puel1} N. {Ben Abdallah}, A. {Mellet}, M. {Puel},
 \emph{Anomalous diffusion limit for kinetic equations with degenerate collision frequency}. Math. Models Methods Appl.  Sci., $21$, $2249$, \oldstylenums{2011}.
 \bibitem{Puel2} N. {Ben Abdallah}, A. {Mellet}, M. {Puel}, \emph{Fractional diffusion limit for collisional kinetic equations: a Hilbert expansion approach}. Kinet. Relat. Models, volume $4$, issue $4$, \oldstylenums{2011}.
  \bibitem{CrousHivertLemou1} N. {Crouseilles}, H. {Hivert}, M. {Lemou}, \emph{Numerical schemes for kinetic equations in the anomalous diffusion limit. {P}art {I}: the 
case of heavy-tailed equilibrium}, arXiv: 1503.04586, \oldstylenums{2015}.
 \bibitem{CrousHivertLemou2} N. {Crouseilles}, H. {Hivert}, M. {Lemou}, \emph{Numerical schemes for kinetic equations in the anomalous diffusion limit. {P}art {II}: the 
case of degenerate collision frequency}, in preparation.
\bibitem{LemouMieussens} M. {Lemou}, L. {Mieussens}, \emph{A new asymptotic preserving scheme based on micro-macro formulation for linear kinetic equations in the diffusion limit}, SIAM J. Sci. Comput, volume $31$ no. $1$, \oldstylenums{2008}.
 \bibitem{MelletMischlerMouhot} A. {Mellet}, S. {Mischler}, C. {Mouhot}, \emph{Fractionnal 
diffusion limit for collisional kinetic equations}, Arch. Ration. Mech. Anal., Vol $199$ 
no.$2$, \oldstylenums{2011}.
\end{thebibliography}
\end{document}